\documentclass[11pt,twoside,a4paper]{article}
\usepackage{amssymb}
\usepackage{amsfonts}
\usepackage{amsmath}
\usepackage{amsthm}
\usepackage{makeidx}
\usepackage[dvips]{graphicx}
\newtheorem{thm}{Theorem}[section]
\newtheorem{dfn}[thm]{Definition}
\newtheorem{lem}[thm]{Lemma}
\newtheorem{cor}[thm]{Corollary}

\newcommand{\F}{\mathcal{F}}
\newcommand{\Ve}{\mathcal{V}}
\newcommand{\Z}{\mathbb{Z}}

\newcommand{\R}{\mathbb{R}}

\newcommand{\conv}{\mathrm{conv}}

\newcommand{\aff}{\mathrm{aff}}

\newcommand{\pair}[2]{\langle #1,#2 \rangle}

\newcommand{\vertsub}{\mathcal{W}}
\newcommand{\Porder}{\mathrm{ord}}
\newcommand{\initsimplex}{I}
\newcommand{\SFP}{\texttt{\large SFP}}
\newcommand{\AddPoint}{\texttt{\large AddPoint}}
\newcommand{\CheckSubset}{\textnormal{\texttt{\large CheckSubset}}}

\def\address#1#2{\begingroup
\noindent\parbox[t]{7.8cm}{%
\small{\scshape\ignorespaces#1}\par\vskip1ex
\noindent\small{\itshape E-mail address}%
\/: #2\par\vskip4ex}\hfill%
\endgroup}%

\title{An algorithm for the classification of smooth Fano polytopes}

\author{Mikkel \O bro}

\begin{document}

\maketitle

\begin{abstract}
We present an algorithm that produces the classification list of smooth Fano $d$-polytopes for any given $d\geq 1$. The input of the algorithm is a single number, namely the positive integer $d$. The algorithm has been used to classify smooth Fano $d$-polytopes for $d\leq 7$. There are 7622 isomorphism classes of smooth Fano $6$-polytopes and 72256 isomorphism classes of smooth Fano $7$-polytopes.
\end{abstract}

\section{Introduction}

Isomorphism classes of smooth toric Fano varieties of dimension $d$
correspond to isomorphism classes of socalled smooth Fano
$d$-polytopes, which are fully dimensional convex lattice polytopes in
$\R^d$, such that the origin is in the interior of the polytopes and
the vertices of every facet is a basis of the integral lattice
$\Z^d\subset \R^d$. Smooth Fano $d$-polytopes have been intensively
studied for the last decades. They have been completely classified up to isomorphism for
$d\leq 4$ (\cite{batyrev82}, \cite{ww82}, \cite{batyrev99}, \cite{sato00}). Under additional assumptions there are classification results valid in every dimension.

To our knowledge smooth Fano $d$-polytopes have been classified in the following cases:
\begin{itemize}
\item When the number of vertices is $d+1$, $d+2$ or $d+3$ (\cite{klein88},\cite{batyrev91}).
\item When the number of vertices is $3d$, which turns out to be the
  upper bound on the number of vertices (\cite{casagrande06}).
\item When the number of vertices is $3d-1$ (\cite{oebro}).
\item When the polytopes are centrally symmetric (\cite{vk85}).
\item When the polytopes are pseudo-symmetric, i.e. there is a facet $F$, such that $-F$ is also a facet (\cite{ewald88}).
\item When there are many pairs of centrally symmetric vertices (\cite{casagrande03}).
\item When the corresponding toric $d$-folds are equipped with an
  extremal contraction, which contracts a toric divisor to a point
  (\cite{bono}) or a curve (\cite{sato03}).
\end{itemize}

Recently a complete classification of smooth Fano $5$-polytopes has been announced (\cite{nillkreuzer}). The approach is to recover smooth Fano $d$-polytopes from their image under the projection along a vertex. This image is a \emph{reflexive} $(d-1)$-polytope (see \cite{batyrev99}), which is a fully-dimensional lattice polytope containing the origin in the interior, such that the dual polytope is also a lattice polytope. Reflexive polytopes have been classified up to dimension 4 using the computer program PALP (\cite{ks98},\cite{ks00}). Using this classification and PALP the authors of \cite{nillkreuzer} succeed in classifying smooth Fano $5$-polytopes.
\\

In this paper we present an algorithm that classifies smooth Fano $d$-polytopes for any given $d\geq 1$. We call this algorithm SFP (for Smooth Fano Polytopes). The input is the positive integer $d$, nothing else is needed. The algorithm has been implemented in C++, and used to classify smooth Fano $d$-polytopes for $d\leq 7$. For $d=6$ and $d=7$ our results are new:

\begin{thm}
There are 7622 isomorphism classes of smooth Fano $6$-polytopes and 72256 isomorphism classes of smooth Fano $7$-polytopes.
\end{thm}

The classification lists of smooth Fano $d$-polytopes, $d\leq 7$, are available on the authors homepage: \texttt{http://home.imf.au.dk/oebro}

A key idea in the algorithm is the notion of a special facet of a smooth Fano $d$-polytope (defined in section \ref{specialfacet}): A facet $F$ of a smooth Fano $d$-polytope is called \emph{special}, if the sum of the vertices of the polytope is a non-negative linear combination of vertices of $F$. This allows us to identify a finite subset $\vertsub_d$ of the lattice $\Z^d$, such that any smooth Fano $d$-polytope is isomorphic to one whose vertices are contained in $\vertsub_d$ (theorem \ref{mainthm}). Thus the problem of classifying smooth Fano $d$-polytopes is reduced to the problem of considering certain subsets of $\vertsub_d$.

We then define a total order on finite subsets of $\Z^d$ and use this to define a total order on the set of smooth Fano $d$-polytopes, which respects isomorphism (section \ref{totalorder}). The SFP-algorithm (described in section \ref{sfpalgo}) goes through certain finite subsets of $\vertsub_d$ in increasing order, and outputs smooth Fano $d$-polytopes in increasing order, such that any smooth Fano $d$-polytope is isomorphic to exactly one in the output list.

As a consequence of the total order on smooth Fano $d$-polytopes, the
algorithm needs not consult the previous output to check for
isomorphism to decide whether or not to output a constructed polytope.

\section{Smooth Fano polytopes}
\label{herearelemmas}

We fix a notation and prove some simple facts about smooth Fano polytopes.

The convex hull of a set $K\in \R^d$ is denoted by $\conv K$. A \emph{polytope} is the convex hull of finitely many points. The dimension of a polytope $P$ is the dimension of the affine hull, $\aff P$, of the polytope $P$. A \emph{$k$-polytope} is a polytope of dimension $k$. A \emph{face} of a polytope is the intersection of a supporting hyperplane with the polytope. Faces of polytopes are polytopes. Faces of dimension 0 are called \emph{vertices}, while faces of codimension 1 and 2 are called \emph{facets} and \emph{ridges}, respectively. The set of vertices of a polytope $P$ is denoted by $\Ve(P)$.

\begin{dfn}
A convex lattice polytope $P$ in $\R^d$ is called a \emph{smooth Fano $d$-polytope}, if the origin is contained in the interior of $P$ and the vertices of every facet of $P$ is a $\Z$-basis of the lattice $\Z^d\subset \R^d$.
\end{dfn}

We consider two smooth Fano $d$-polytopes $P_1,P_2$ to be \emph{isomorphic}, if there exists a bijective linear map $\varphi:\R^d\to \R^d$, such that $\varphi(\Z^d)=\Z^d$ and $\varphi(P_1)=P_2$.

Whenever $F$ is a $(d-1)$-simplex in $\R^d$, such that $0\notin \aff F$, we let $u_F\in (\R^d)^*$ be the unique element determined by $\pair{u_F}{F}=\{1\}$. For every $w\in \Ve(F)$ we define $u_F^w\in (\R^d)^*$ to be the element where $\pair{u_F^w}{w}=1$ and $\pair{u_F^w}{w'}=0$ for every $w'\in \Ve(F)$, $w'\neq w$. Then $\{u_F^w|w\in\Ve(F)\}$ is the basis of $(\R^d)^*$ dual to the basis $\Ve(F)$ of $\R^d$.

When $F$ is a facet of a smooth Fano polytope and $v\in\Ve(P)$, we certainly have $\pair{u_F}{v}\in\Z$ and
$$
\pair{u_F}{v}=1\ \ \Longleftrightarrow \ \ v\in\Ve(F)\ \ \ \textnormal{ and }\ \ \ \pair{u_F}{v}\leq 0\ \ \Longleftrightarrow \ \ v\notin \Ve(F).
$$ 
The lemma below concerns the relation between the elements $u_F$ and $u_{F'}$, when $F$ and $F'$ are adjacent facets.

\begin{lem}
\label{firstlemma}
Let $F$ be a facet of a smooth Fano polytope $P$ and $v\in \Ve(F)$. Let $F'$ be the unique facet which intersects $F$ in a ridge $R$ of $P$, $v\notin \Ve(R)$. Let $v'=\Ve(F')\setminus \Ve(R)$.

Then
\begin{enumerate}
\item $\pair{u_F^v}{v'}=-1$.
\label{firstlemma1}
\item $\pair{u_F}{v'}=\pair{u_{F'}}{v}$.
\item $\pair{u_{F'}}{x}=\pair{u_F}{x}+\pair{u_F^v}{x}(\pair{u_F}{v'}-1)$ for any $x\in\R^d$.
\item In particular,
\begin{itemize}
\item $\pair{u_F^v}{x}<0$ iff $\pair{u_{F'}}{x}>\pair{u_F}{x}$.
\item $\pair{u_F^v}{x}>0$ iff $\pair{u_{F'}}{x}<\pair{u_F}{x}$.
\item $\pair{u_F^v}{x}=0$ iff $\pair{u_{F'}}{x}=\pair{u_F}{x}$.
\end{itemize}
for any $x\in\R^d$.
\item Suppose $x\neq v'$ is a vertex of $P$ where $\pair{u_F^v}{x}<0$. Then $\pair{u_F}{v'}>\pair{u_F}{x}$.
\label{closeneighbor}
\end{enumerate}
\begin{proof}
The sets $\Ve(F)$ and $\Ve(F')$ are both bases of the lattice $\Z^d$ and the first statement follows.

We have $v+v'\in\textnormal{span} (F\cap F')$, and then the second statement follows.

Use the previous statements to calculate $\pair{u_{F'}}{x}$.
\begin{eqnarray*}
\pair{u_{F'}}{x}&=& \pair{u_{F'}}{\sum_{w\in\Ve(F)} \pair{u_F^w}{x} w}\\
&=& \sum_{w\in\Ve(F)\setminus\{v\}} \pair{u_F^w}{x} + \pair{u_F^v}{x} \pair{u_{F'}}{v}\\
&=& \pair{u_F}{x} + \pair{u_F^v}{x} \big( \pair{u_{F'}}{v}-1\big)\\
&=& \pair{u_F}{x}+ \pair{u_F^v}{x} \big( \pair{u_{F}}{v'} -1 \big).
\end{eqnarray*}
As $\pair{u_F}{v'}-1<0$ the three equivalences follow directly.

Suppose there is a vertex $x\in\Ve(P)$, such that $\pair{u_F^v}{x}<0$ and $\pair{u_F}{v'}\leq \pair{u_F}{x}$. Then
$$
\pair{u_{F'}}{x}=\pair{u_F}{x}+\pair{u_F^v}{x}(\pair{u_F}{v'}-1)\geq \pair{u_F}{x}-(\pair{u_F}{v'}-1)\geq 1.
$$
Hence $x$ is on the facet $F'$. But this cannot be the case as $\Ve(F')=\{v'\}\cup \Ve(F)\setminus \{v\}$. Thus no such $x$ exists.

And we're done.
\end{proof}
\end{lem}

In the next lemma we show a lower bound on the numbers $\pair{u_F^w}{v}$, $w\in \Ve(F)$, for any facet $F$ and any vertex $v$ of a smooth Fano $d$-polytope.

\begin{lem} Let $F$ be a facet and $v$ a vertex of a smooth Fano polytope $P$. Then
$$
\pair{u_F^w}{v}\geq \left\{ \begin{array}{cl}
0 & \pair{u_F}{v}=1\\
-1 & \pair{u_F}{v}=0\\
\pair{u_F}{v} & \pair{u_F}{v}<0 \end{array} \right.
$$
for every $w\in\Ve(F)$.
\label{coef_lemma}
\begin{proof}
When $\pair{u_F}{v}=1$ the statement is obvious.

Suppose $\pair{u_F}{v}=0$ and $\pair{u_F^w}{v}<0$ for some $w\in\Ve(F)$. Let $F'$ be the unique facet intersecting $F$ in the ridge $\conv\{ \Ve(F)\setminus \{w\}\}$. By lemma \ref{firstlemma} $\pair{u_{F'}}{v}>0$. As $\pair{u_{F'}}{v}\in \Z$ we must have $\pair{u_{F'}}{v}=1$. This implies $\pair{u_F}{v}=-1$.

Suppose $\pair{u_F}{v}<0$ and $\pair{u_F^w}{v}<\pair{u_F}{v}\leq -1$ for some $w\in\Ve(F)$. Let $F'\neq F$ be the facet containing the ridge $\conv\{ \Ve(F)\setminus \{w\}\}$, and let $w'$ be the unique vertex in $\Ve(F')\setminus \Ve(F)$. Then by lemma \ref{firstlemma}
$$
\pair{u_{F'}}{v}=\pair{u_F}{v}+\pair{u_F^w}{v}(\pair{u_F}{w'}-1)\geq \pair{u_F}{v}-\pair{u_F^w}{v}.
$$
If $\pair{u_F}{v}-\pair{u_F^w}{v}>0$, then $v$ is on the facet $F'$. But this is not the case as $\pair{u_F^w}{v}<-1$. We conclude that $\pair{u_F^w}{v}\geq \pair{u_F}{v}$.
\end{proof}
\end{lem}

When $F$ is a facet and $v$ a vertex of a smooth Fano $d$-polytope $P$, such that $\pair{u_F}{v}=0$, we can say something about the face lattice of $P$. 

\begin{lem}[\cite{debarre} section 2.3 remark 5(2), \cite{nill05} lemma 5.5]
\label{hyper0lemma}
Let $F$ be a facet and $v$ be
  vertex of a smooth Fano polytope $P$. Suppose $\pair{u_F}{v}=0$.

Then $\conv\{ \{v\}\cup\Ve(F)\setminus\{w\}\}$ is a facet of $P$ for every $w\in\Ve(F)$ with $\pair{u_F^w}{v}=-1$.
\begin{proof}
Follows from the proof of lemma \ref{coef_lemma}.
\end{proof}
\end{lem}

\section{Special embeddings of smooth Fano polytopes}

In this section we find a concrete finite subset $\vertsub_d$ of $\Z^d$ with the nice property that any smooth Fano $d$-polytope is isomorphic to one whose vertices are contained in $\vertsub_d$. The problem of classifying smooth Fano $d$-polytopes is then reduced to considering subsets of $\vertsub_d$.

\subsection{Special facets}

\label{specialfacet}
The following definition is a key concept.
\begin{dfn}
A facet $F$ of a smooth Fano $d$-polytope $P$ is called \emph{special}, if the sum of the vertices of $P$ is a non-negative linear combination of $\Ve(F)$, that is
$$
\sum_{v\in\Ve(P)} v=\sum_{w\in\Ve(F)} a_w w\ \ ,\ a_w\geq 0.
$$
\end{dfn}
Clearly, any smooth Fano $d$-polytope has at least one special facet.

Let $F$ be a special facet of a smooth Fano $d$-polytope $P$. Then
$$
0\leq \pair{u_F}{\sum_{v\in\Ve(P)} v} = d+\sum_{v\in\Ve(P),\pair{u_F}{v}<0} \pair{u_F}{v},
$$
which implies $-d\leq \pair{u_F}{v}\leq 1$ for any vertex $v$ of $P$. By using the lower bound on the numbers $\pair{u_F^w}{v}$, $w\in\Ve(F)$ (see lemma \ref{coef_lemma}), we can find an explicite finite subset of the lattice $\Z^d$, such that every $v\in\Ve(P)$ is contained in this subset. In the following lemma we generalize this observation to subsets of $\Ve(P)$ containing $\Ve(F)$.

\begin{lem}
Let $P$ be a smooth Fano polytope. Let $F$ be a special facet of $P$ and let $V$ be a subset of $\Ve (P)$ containing $\Ve(F)$, whose sum is $\nu$.
$$
\nu=\sum_{v\in V} v.
$$
Then
$$
\pair{u_F}{\nu} \geq 0
$$
and
$$
\pair{u_F^w}{\nu} \leq \pair{u_F}{\nu}+1
$$
for every $w\in \Ve(F)$.
\label{sumlem}
\begin{proof}
For convenience we set $U=\Ve(P)\setminus V$ and $\mu=\sum_{v\in U} v$. Since $F$ is a special facet we know that
$$
0\leq \pair{u_F}{\sum_{v\in \Ve(P)} v}=\pair{u_F}{\nu}+\pair{u_F}{\mu}.
$$
The set $\Ve(F)$ is contained in $V$ so $\pair{u_F}{v}\leq 0$ for every $v$ in $U$, hence $\pair{u_F}{\nu}\geq 0$.

Suppose that for some $w\in \Ve (F)$ we have
$\pair{u_F^w}{\nu}>\pair{u_F}{\nu}+1$. By lemma \ref{coef_lemma} we know that
$$
\pair{u_F^w}{v}\geq \left\{ \begin{array}{ll}
-1 & \pair{u_F}{v}=0\\
\pair{u_F}{v} & \pair{u_F}{v}<0
\end{array}
\right.
$$
for every vertex $v\in\Ve(P)\setminus \Ve(F)$. There is at most one vertex $v$ of $P$, $\pair{u_F}{v}=0$, with negative coefficient $\pair{u_F^w}{v}$ (lemma \ref{hyper0lemma}). So
$$
\pair{u_F^w}{\mu}\geq \pair{u_F}{\mu}-1.
$$
Now, consider $\pair{u_F^w}{\sum_{v\in\Ve(P)} v}$.
$$
\pair{u_F^w}{\sum_{v\in\Ve(P)} v}=\pair{u_F^w}{\nu}+\pair{u_F^w}{\mu}>\pair{u_F}{\nu}+\pair{u_F}{\mu}=\pair{u_F}{\sum_{v\in\Ve(P)} v}.
$$
But this implies that $\pair{u_F^x}{\sum_{v\in\Ve(P)} v}$ is negative for some $x\in\Ve(F)$. A contradiction.
\end{proof}
\end{lem}

\begin{cor}
\label{coefcor}
Let $F$ be a special facet and $v$ any vertex of a smooth Fano $d$-polytope. Then $-d\leq \pair{u_F}{v}\leq 1$ and 
$$
\left.
\begin{array}{c}
0 \\
-1\\
\pair{u_F}{v}
\end{array}
\right\}
\leq \pair{u_F^w}{v} \leq \left\{
\begin{array}{c @{\ \ ,\ } l}
1 & \pair{u_F}{v}=1\\
d-1 & \pair{u_F}{v}=0\\
d+\pair{u_F}{v} & \pair{u_F}{v}<0\\
\end{array} \right.
$$
for every $w\in\Ve(F)$.
\begin{proof}
For $\pair{u_F}{v}=1$ the statement is obvious. When $\pair{u_F}{v}=0$ the coefficients of $v$ with respect to the basis $\Ve(F)$ is bounded below by $-1$ (lemma \ref{coef_lemma}), so no coefficient exceeds $d-1$.

So the case $\pair{u_F}{v}<0$ remains. The lower bound is by lemma
\ref{coef_lemma}. Use lemma \ref{sumlem} on the subset $V=\Ve(F)\cup
\{v\}$ to prove the upper bound.
\end{proof}
\end{cor}

\subsection{Special embeddings}

Let $(e_1,\ldots,e_d)$ be a fixed basis of the lattice $\Z^d\subset \R^d$.

\begin{dfn}
Let $P$ be a smooth Fano $d$-polytope. Any smooth Fano $d$-polytope $Q$, with $\conv\{e_1,\ldots,e_d\}$ as a special facet, is called a \emph{special embedding of $P$}, if $P$ and $Q$ are isomorphic.
\end{dfn}

Obviously, for any smooth Fano polytope $P$, there exists at least one special embedding of $P$. As any polytope has finitely many facets, there exists only finitely many special embeddings of $P$.

Now we define a subset of $\Z^d$ which will play an important part in what follows.
\begin{dfn}
By $\vertsub_d$ we denote the maximal set (with respect to inclusion) of lattice points in $\Z^d$ such that
\begin{enumerate}
\item The origin is not contained in $\vertsub_d$.
\item The points in $\vertsub_d$ are primitive lattice points.
\item If $a_1e_1+\ldots+a_de_d\in \vertsub_d$, then $-d\leq a\leq 1$ for $a=a_1+\ldots+a_d$ and
$$
\left.
\begin{array}{c}
0 \\
-1\\
a
\end{array}
\right\}
\leq a_i \leq \left\{
\begin{array}{c @{\ \ ,\ } l}
1 & a=1\\
d-1 & a=0\\
d+a & a<0\\
\end{array} \right.
$$
for every $i=1,\ldots,d$.
\end{enumerate}
\end{dfn}

The next theorem is one of the key results in this paper. It allows us to classify smooth Fano $d$-polytopes by considering subsets of the explicitely given set $\vertsub_d$.

\begin{thm}
\label{mainthm}
Let $P$ be an arbitrary smooth Fano $d$-polytope, and $Q$ any special embedding of $P$. Then $\Ve(Q)$ is contained in the set $\vertsub_d$.
\begin{proof}
Follows directly from corollary \ref{coefcor} and the definition of $\vertsub_d$.
\end{proof}
\end{thm}

\section{Total ordering of smooth Fano polytopes}
\label{totalorder}
In this section we define a total order on the set of smooth Fano $d$-polytopes for any fixed $d\geq 1$.

Throughout the section $(e_1,\ldots,e_d)$ is a fixed basis of the lattice $\Z^d$.

\subsection{The order of a lattice point}

We begin by defining a total order $\preceq$ on $\Z^d$.
\begin{dfn}
Let $x=x_1e_1+\ldots+x_de_d,\ y=y_1e_1+\ldots+y_de_d$ be two lattice points in $\Z^d$. We define $x\preceq y$ if and only if
$$
(-x_1-\ldots-x_d,x_1,\ldots,x_d)\leq_{lex} (-y_1-\ldots-y_d,y_1,\ldots,y_d),
$$
where $\leq_{lex}$ is the lexicographical ordering on the product of $d+1$ copies of the ordered set $(\Z,\leq)$.

The ordering $\preceq$ is a total order on $\Z^d$.

\end{dfn}

\textbf{Example.} $(0,1)\prec (-1,1) \prec (1,-1)\prec (-1,0)$.
\\

Let $V$ be any nonempty finite subset of lattice points in $\Z^d$. We define $\max V$ to the maximal element in $V$ with respect to the ordering $\preceq$. Similarly, $\min V$ is defined to be the minimal element in $V$.

A important property of the ordering is shown in the following lemma.
\begin{lem}
\label{proporder}
Let $P$ be a smooth Fano $d$-polytope, such that $\conv\{e_1,\ldots, e_d\}$ is a facet of $P$. For every $1\leq i\leq d$, let $v_i\neq e_i$ denote the vertex of $P$, such that $\conv\{e_1,\ldots,e_{i-1},v_i,e_{i+1},\ldots,e_d\}$ is a facet of $P$.

Then $v_i=\min \{v\in\Ve(P)\ |\ \pair{u_F^{e_i}}{v}<0\}$.
\begin{proof}
By lemma \ref{firstlemma}.(\ref{firstlemma1}) the vertex $v_i$ is in the set $\{v\in\Ve(P)\ |\ \pair{u_F^{e_i}}{v}<0\}$, and by lemma \ref{firstlemma}.(\ref{closeneighbor}) and the definition of the ordering $\preceq$, $v_i$ is the minimal element in this set.
\end{proof}
\end{lem}
In fact, we have chosen the ordering $\preceq$ to obtain the property of lemma \ref{proporder}, and any other total order on $\Z^d$ having this property can be used in what follows.

\subsection{The order of a smooth Fano $d$-polytope}

We can now define an ordering on finite subsets of $\Z^d$. The ordering is defined recursively.
\begin{dfn}
Let $X$ and $Y$ be finite subsets of $\Z^d$. We define $X\preceq Y$ if and only if $X=\emptyset$ or
$$
Y\neq \emptyset \ \land \ (\min X\prec \min Y \ \lor \ (\min X=\min Y\ \land X\setminus\{\min X\}\preceq Y\setminus \{\min Y\}) ).
$$
\end{dfn}

\textbf{Example. } $\emptyset\prec \{(0,1)\}\prec \{(0,1),(-1,1)\}\prec \{(0,1),(1,-1)\}\prec \{(-1,1)\}.$
\\

When $W$ is a nonempty finite set of subsets of $\Z^d$, we define $\max W$ to be the maximal element in $W$ with respect to the ordering of subsets $\preceq$. Similarly, $\min W$ is the minimal element in $W$.

Now, we are ready to define the order of a smooth Fano $d$-polytope.

\begin{dfn}
Let $P$ be a smooth Fano $d$-polytope. The order of $P$, $\Porder(P)$, is defined as
$$
\Porder(P):=\min \{ \Ve(Q)\ |\ Q\ \textnormal{a special embedding of $P$}\}.
$$
The set is non-empty and finite, so $\Porder(P)$ is well-defined.

Let $P_1$ and $P_2$ be two smooth Fano $d$-polytopes. We say that $P_1\leq P_2$ if and only if $\Porder(P_1)\preceq \Porder(P_2)$. This is indeed a total order on the set of isomorphism classes of smooth Fano $d$-polytopes.
\end{dfn}

\subsection{Permutation of basisvectors and presubsets}

The group $S_d$ of permutations of $d$ elements acts on $\Z^d$ is the obvious way by permuting the basisvectors:
$$
\sigma.(a_1e_1+\ldots+a_de_d):=a_1e_{\sigma(1)}+\ldots+a_de_{\sigma(d)}\ \ ,\ \sigma\in S_d.
$$
Similarly, $S_d$ acts on subsets of $\Z^d$:
$$
\sigma.X:=\{\sigma.x\ |\ x\in X\}.
$$
In this notation we clearly have for any special embedding $P$ of a smooth Fano $d$-polytope
$$
\Porder(P)\preceq \min \{\sigma.\Ve(P)\ |\ \sigma\in S_d\}.
$$
Let $V$ and $W$ be finite subsets of $\Z^d$. We say that $V$ is a \emph{presubset} of $W$, if $V\subseteq W$ and $v\prec w$ whenever $v\in V$ and $w\in W\setminus V$.
\\

\textbf{Example.} $\{(0,1),(-1,1)\}$ is a presubset of $\{(0,1),(-1,1),(1,-1)\}$, while $\{(0,1),(1,-1)\}$ is not.

\begin{lem}
Let $P$ be a smooth Fano polytope. Then every presubset $V$ of $\Porder(P)$ is the minimal element in $\{\sigma.V\ |\ \sigma\in S_d\}$.
\label{permutepresubset}
\begin{proof}
Let $\Porder(P)=\{v_1,\ldots,v_n\}$, $v_1\prec \ldots \prec v_n$. Suppose there exists a permutation $\sigma$ and a $k$, $1\leq k\leq n$, such that
$$
\sigma.\{v_1,\ldots,v_k\}=\{w_1,\ldots,w_k\}\prec \{v_1,\ldots,v_k\},
$$
where $w_1\prec \ldots \prec w_k$. Then there is a number $j$, $1\leq j\leq k$, such that $w_i=v_i$ for every $1\leq i<j$ and $w_j\prec v_j$.

Let $\sigma$ act on $\{v_1,\ldots,v_n\}$.
$$
\sigma.\{v_1,\ldots,v_n\}=\{x_1,\ldots,x_n\}\ \ ,\ x_1\prec \ldots \prec x_n.
$$
Then $x_i\preceq v_i$ for every $1\leq i< j$ and $x_j\prec v_j$. So $\sigma.\Porder(P)\prec \Porder(P)$, but this contradicts the definition of $\Porder(P)$.
\end{proof}
\end{lem}

\section{The SFP-algorithm}
\label{sfpalgo}

In this section we describe an algorithm that produces the classification list of smooth Fano $d$-polytopes for any given $d\geq 1$. The algorithm works by going through certain finite subsets of $\vertsub_d$ in increasing order (with respect to the ordering defined in the previous section). It will output a subset $V$ iff $\conv V$ is a smooth Fano $d$-polytope $P$ and $\Porder(P)=V$.

Throughout the whole section $(e_1,\ldots,e_d)$ is a fixed basis of $\Z^d$ and $\initsimplex$ denotes the $(d-1)$-simplex $\conv\{e_1,\ldots,e_d\}$.

\subsection{The SFP-algorithm}

The SFP-algorithm consists of three functions,
\begin{center}
\SFP, \AddPoint \ and \CheckSubset.
\end{center}
The finite subsets of $\vertsub_d$ are constructed by the function
\AddPoint , which takes a subset $V$, $\{e_1,\ldots,e_d\}\subseteq
V\subseteq \vertsub_d$, together with a finite set $\F$, $I\in\F$, of $(d-1)$-simplices in $\R^d$ as input. It then goes through every $v$ in the set
$$
\{v\in\vertsub_d\ |\ \max V\prec v\}
$$
in increasing order, and recursively calls itself with input $V\cup \{v\}$ and some set $\F'$ of $(d-1)$-simplices of $\R^d$, $\F\subseteq \F'$. In this way subsets of $\vertsub_d$ are considered in increasing order.

Whenever \AddPoint \ is called, it checks if the input set $V$ is the vertex set of a special embedding of a smooth Fano $d$-polytope $P$ such that $\Porder(P)=V$, in which case the polytope $P=\conv V$ is outputted.

For any given integer $d\geq 1$ the function \SFP \ calls the function \AddPoint \ with input $\{e_1,\ldots,e_d\}$ and $\{\initsimplex\}$. In this way a call \SFP$(d)$ will make the algorithm go through every finite subset of $\vertsub_d$ containing $\{e_1,\ldots,e_d\}$, and smooth Fano $d$-polytopes are outputted in strictly increasing order.

It is vital for the effectiveness of the SFP-algorithm, that there is some efficient way to check if a subset $V\subseteq \vertsub_d$ is a presubset of $\Porder (P)$ for some smooth Fano $d$-polytope $P$. The function \AddPoint \ should perform this check before the recursive call \AddPoint $(V,\F')$.

If $P$ is any smooth Fano $d$-polytope, then any presubset $V$ of $\Porder(P)$ is the minimal element in the set $\{\sigma.V|\sigma\in S_d\}$ (by lemma \ref{permutepresubset}). In other words, if there exists a permutation $\sigma$ such that $\sigma.V\prec V$, then the algorithm should not make the recursive call \AddPoint $(V)$.

But this is not the only test we wish to perform on a subset $V$
before the recursive call. The function \CheckSubset \ performs
another test: It takes a subset $V$, $\{e_1,\ldots,e_d\}\subseteq
V\subseteq \vertsub_d$ as input together with a finite set of
$(d-1)$-simplices $\F$, $I\in \F$, and returns a set $\F'$ of $(d-1)$-simplices containing $\F$, if there exists a special embedding $P$ of a smooth Fano $d$-polytope, such that
\begin{enumerate}
\item $V$ is a presubset of $\Ve(P)$
\item $\F$ is a subset of the facets of $P$
\end{enumerate}
This is proved in theorem \ref{justify1}. If no such special embedding exists, then \CheckSubset \ returns false in many cases, but not always! Only when \CheckSubset$(V,\F)$ returns a set $\F'$ of simplices, we allow the recursive call \AddPoint $(V,\F')$.

Given input $V\subseteq \vertsub_d$ and a set $\F$ of $(d-1)$-simplices of $\R^d$, the function \CheckSubset \ works in the following way: Suppose $V$ is a presubset of $\Ve(P)$ for some special embedding $P$ of a smooth Fano $d$-polytope and $\F$ is a subset of the facets of $P$. Deduce as much as possible of the face lattice of $P$ and look for contradictions to the lemmas stated in section \ref{herearelemmas}. The more facets we know of $P$, the more restrictions we can put on the vertex set $\Ve(P)$, and then on $V$. If a contradiction arises, return false. Otherwise, return the deduced set of facets of $P$.

The following example illustrates how the function \CheckSubset \ works.

\subsection{An example of the reasoning in \CheckSubset}

Let $d=5$ and $V=\{v_1,\ldots,v_8\}$, where
\begin{center}
$v_1=e_1$ , $v_2=e_2$ , $v_3=e_3$ , $v_4=e_4$ , $v_5=e_5$

$v_6=-e_1-e_2+e_4+e_5$ , $v_7=e_2-e_3-e_4$ , $v_8=-e_4-e_5$.
\end{center}
Suppose $P$ is a special embedding of a smooth Fano $5$-polytope, such that $V$ is a presubset of $\Ve(P)$. Certainly, the simplex $\initsimplex$ is a facet of $P$. 

Notice, that $V$ does not violate lemma \ref{sumlem}.
$$
v_1+\ldots +v_8=e_2+e_5.
$$
If $V$ did contradict lemma \ref{sumlem}, then the polytope $P$ could not exist, and \CheckSubset$(V,\{\initsimplex\})$ should return false.

For simplicity we denote any $k$-simplex $\conv\{v_{i_1},\ldots,v_{i_k}\}$ by $\{i_1,\ldots,i_k\}$. 

Since $\pair{u_\initsimplex}{v_6}=0$, the simplices $F_1=\{2,3,4,5,6\}$ and $F_2=\{1,3,4,5,6\}$ are facets of $P$ (lemma \ref{hyper0lemma}).

There are exactly two facets of $P$ containing the ridge $\{1,2,4,5\}$. One of them is $\initsimplex$. Suppose the other one is $\{1,2,4,5,9\}$, where $v_9$ is some lattice point not in $V$, $v_9\in\Ve(P)$. Then $\pair{u_{\initsimplex}}{v_9}>\pair{u_{\initsimplex}}{v_7}$ by lemma \ref{firstlemma}.(\ref{closeneighbor}) and then $v_9\prec v_7$ by the definition of the ordering of lattice points $\Z^d$. But then $V$ is not a presubset of $\Ve(P)$. This is the nice property of the ordering of $\Z^d$, and the reason why we chose it as we did. We conclude that $F_3=\{1,2,4,5,7\}$ is a facet of $P$, and by similar reasoning $F_4=\{1,2,3,5,8\}$ and $F_5=\{1,2,3,4,8\}$ are facets of $P$. 

Now, for each of the facets $F_i$ and every point $v_j\in V$, we check if $\pair{u_{F_i}}{v_j}=0$. If this is the case, then by lemma \ref{hyper0lemma} $\conv(\{v_j\}\cup \Ve(F_i)\setminus\{w\})$ is a facet of $P$ for every $w\in\Ve(F_i)$ where $\pair{u_{F_i}^w}{v_j}<0$. In this way we get that
$$
\{2,4,5,6,7\}\ ,\ \{1,4,5,6,7\}\ ,\ \{1,2,3,7,8\}\ ,\ \{1,3,5,7,8\}
$$
are facets of $P$.

We continue in this way, until we cannot deduce any new facet of $P$. Every time we find a new facet $F$ we check that $v$ is beneath $F$ (that is $\pair{u_F}{v}\leq 1$) and that lemma \ref{coef_lemma} holds for any $v\in V$. If not, then \CheckSubset$(V,\{\initsimplex\})$ should return false.

If no contradiction arises, \CheckSubset$(V,\{\initsimplex\})$ returns the set of deduced facets.

\subsection{The SFP-algorithm in pseudo-code}

Input: A positive integer $d$.

Output: A list of special embeddings of smooth Fano $d$-polytopes, such that
\begin{enumerate}
\item Any smooth Fano $d$-polytope is isomorphic to one and only one polytope in the output list.
\item If $P$ is a smooth Fano $d$-polytope in the output list, then $\Ve(P)=\Porder(P)$.
\item If $P_1$ and $P_2$ are two non-isomorphic smooth Fano $d$-polytopes in the output list and $P_1$ preceeds $P_2$ in the output list, then $\Porder(P_1)\prec \Porder(P_2)$.
\end{enumerate}

\SFP \ ( an integer $d\geq 1$ )
\begin{enumerate}
\item Construct the set $V=\{e_1,\ldots,e_d\}$ and the simplex $\initsimplex=\conv V$.
\item Call the function \AddPoint $(V,\{\initsimplex\})$.
\item End program.
\end{enumerate}

\AddPoint \ ( a subset $V$ where $\{e_1,\ldots,e_d\}\subseteq V\subseteq
\vertsub_d$ , a set of $(d-1)$-simplices $\F$ in $\R^d$ where $I\in\F$
)
\begin{enumerate}
\item If $P=\conv (\Ve(V))$ is a smooth Fano $d$-polytope and $\Ve(V)=\Porder(P)$, then output $P$.
\label{outputstep}
\item Go through every $v\in\vertsub_d$, $\max \Ve(V)\prec v$, in increasing order with respect to the ordering $\prec$:
\label{theloop}
\begin{enumerate}
\item If \CheckSubset$(V\cup \{v\},\F)$ returns false, then goto (d). Otherwise let $\F'$ be the returned set of $(d-1)$-simplices.
\label{addtest1}
\item If $V\cup\{v\}\neq \min \{\sigma.(V\cup\{v\})\ |\ \sigma\in S_d\}$, then goto (d).
\label{addtest2}
\item Call the function \AddPoint ($V\cup\{v\},\F'$).
\item Let $v$ be the next element in $\vertsub_d$ and go back to (a).
\end{enumerate}
\item Return
\end{enumerate}

\CheckSubset \ ( a subset $V$ where $\{e_1,\ldots,e_d\}\subseteq V\subseteq
\vertsub_d$ , a set of $(d-1)$-simplices $\F$ in $\R^d$ where $I\in\F$
)
\begin{enumerate}
\item Let $\nu=\sum_{v\in V} v$.
\item If $\pair{u_I}{\nu}<0$, then return false.
\label{test1}
\item If $\pair{u_I^{e_i}}{\nu}>1+\pair{u_I}{\nu}$ for some $i$, then return false.
\label{test2}
\item Let $\F'=\F$.
\item For every $i\in\{1,\ldots,d\}$: If the set $\{v\in V|\pair{u_{\initsimplex}^{e_i}}{v}<0\}$ is equal to $\{\max V\}$, then add the simplex $\conv(\{\max V\}\cup \Ve(\initsimplex)\setminus\{e_i\})$ to $\F'$.
\label{addneighbors}
\item If there exists $F\in\F'$ such that $\Ve(F)$ is not a $\Z$-basis of $\Z^d$, then return false.
\label{loopbegin}
\item If there exists $F\in\F'$ and $v\in V$ such that $\pair{u_F}{v}>1$, then return false.
\item If there exists $F\in\F'$, $v\in V$ and $w\in\Ve(F)$, such that 
$$
\pair{u_F^w}{v}<\left\{ \begin{array}{cl}
0 & \pair{u_F}{v}=1\\
-1 & \pair{u_F}{v}=0\\
\pair{u_F}{v} & \pair{u_F}{v}<0 \end{array} \right.
$$
then return false.
\label{lasttest}
\item If there exists $F\in\F'$, $v\in V$ and $w\in\Ve(F)$, such that $\pair{u_F}{v}=0$ and $\pair{u_F^w}{v}=-1$, then consider the simplex $F'=\conv(\{v\}\cup \Ve(F)\setminus \{w\})$. If $F'\notin \F'$, then add $F'$ to $\F'$ and go back to step \ref{loopbegin}.
\label{addnewsimplices}
\item Return $\F'$.
\end{enumerate}

\subsection{Justification of the SFP-algorithm}

The following theorems justify the SFP-algorithm.

\begin{thm}
\label{justify1}
Let $P$ be a special embedding of a smooth Fano $d$-polytope and $V$ a presubset of $\Ve(P)$, such that $\{e_1,\ldots,e_d\}\subseteq V$. Let $\F$ be a set of facets of $P$.

Then \CheckSubset$(V,\F)$ returns a subset $\F'$ of the facets of $P$ and $\F\subseteq \F'$.
\begin{proof}
By lemma \ref{sumlem} the subset $V$ will pass the tests in step
\ref{test1} and \ref{test2} in \CheckSubset .

The function \CheckSubset \ constructs a set $\F'$ of $(d-1)$-simplices containing the input set $\F$. We now wish to prove that every simplex $F$ in $\F'$ is a facet of $P$: By the assumptions the subset $\F\subseteq \F'$ consists of facets of $P$.

Consider the addition of a simplex $F_i$, $1\leq i\leq d$, in step \ref{addneighbors}:
$$
F_i=\conv(\{\max V\} \cup \Ve(\initsimplex)\setminus \{e_i\}).
$$
As $\max V$ is the only element in the set $\{v\in V|\pair{u_{\initsimplex}^{e_i}}{v}<0\}$ and $V$ is a presubset of $\Ve(P)$, $F_i$ is a facet of $P$ by lemma \ref{proporder}.

Consider the addition of simplices in step \ref{addnewsimplices}: If $F$ is a facet of $P$, then by lemma \ref{hyper0lemma} the simplex $\conv(\{v\}\cup \Ve(F)\setminus \{w\})$ is a facet of $P$.

By induction we conclude, that every simplex in $\F'$ is a facet of $P$. Then any simplex $F\in \F'$ will pass the tests in steps \ref{loopbegin}--\ref{lasttest} (use lemma \ref{coef_lemma} to see that the last test is passed).

This proves the theorem.
\end{proof}
\end{thm}

\begin{thm}
\label{justify2}
The SFP-algorithm produces the promised output.
\begin{proof}
Let $P$ be a smooth Fano $d$-polytope. Clearly, $P$ is isomorphic to at most one polytope in the output list.

Let $Q$ be a special embedding of $P$ such that $\Ve(Q)=\Porder(P)$. We need to show that $Q$ is in the output list. Let $\Ve(Q)=\{e_1,\ldots,e_d,q_1,\ldots,q_k\}$, where $q_1\prec \ldots \prec q_k$, and let $V_i=\{e_1,\ldots,e_d,q_1,\ldots,q_i\}$ for every $1\leq i\leq k$.

Certainly the function \AddPoint \ has been called with input $\{e_1,\ldots,e_d\}$ and $\{\initsimplex\}$.

By theorem \ref{justify1} the function call \CheckSubset$(V_1,\{\initsimplex\})$ returns a set $\F_1$ of $(d-1)$-simplices which are facets of $Q$, $\initsimplex\subset \F_1$. By lemma \ref{permutepresubset} the set $V_1$ passes the test in \ref{addtest2} in \AddPoint . Then \AddPoint \ is called recursively with input $V_1$ and $\F_1$.

The call \CheckSubset$(V_1,\F_1)$ returns a subset $\F_2$ of facets of $Q$, and the set $V_2$ passes the test in \ref{addtest2} in \AddPoint . So the call \AddPoint $(V_2,\F_2)$ is made.

Proceed in this way to see that the call \AddPoint $(V_k,\F_k)$ is made, and then the polytope $Q=\conv V_k$ is outputted in step \ref{outputstep} in \AddPoint .
\end{proof}
\end{thm}

\section{Classification results and where to get them}

A modified version of the SFP-algorithm has been implemented in C++, and used to classify smooth Fano $d$-polytopes for $d\leq 7$. On an average home computer our program needs less than one day (january 2007) to construct the classification list of smooth Fano $7$-polytopes. These lists can be downloaded from the authors homepage: \texttt{http://home.imf.au.dk/oebro}

An advantage of the SFP-algorithm is that it requires almost no memory: When the algorithm has found a smooth Fano $d$-polytope $P$, it needs not consult the output list to decide whether to output the polytope $P$ or not. The construction guarentees that $\Ve(P)=\min\{\sigma.\Ve(P)\ |\ \sigma\in S_d\}$ and it remains to check if $\Ve(P)=\Porder(P)$. Thus there is no need of storing the output list.

The table below shows the number of isomorphism classes of smooth Fano $d$-polytopes with $n$ vertices.

\begin{center}
\begin{tabular}{|c|c|c|c|c|c|c|c|}
\hline
$n$ & $d=1$ & $d=2$ & $d=3$ & $d=4$ & $d=5$ & $d=6$ & $d=7$\\
\hline
1 & & & & & & &\\
2 & 1 & & & & & &\\
3 & & 1 & & & & &\\
4 & & 2 & 1 & & & &\\
5 & & 1 & 4 & 1 & & &\\
6 & & 1 & 7 & 9 & 1 & &\\
7 & &   & 4 & 28 & 15 & 1 &\\
8 & &   & 2 & 47 & 91 & 26 & 1\\
9 & &   &   & 27 & 268& 257 & 40\\
10& &   &   & 10 & 312& 1318 & 643\\
11& &   &   &  1 & 137& 2807 & 5347\\
12& &   &   &  1 & 35 & 2204 & 19516\\
13& &   &   &    & 5  & 771 & 26312\\
14& &   &   &    & 2  & 186 & 14758\\
15& &   &   &    &    & 39 & 4362\\
16& &   &   &    &    & 11 & 1013\\
17& &   &   &    &    & 1 & 214\\
18& &   &   &    &    & 1 & 43\\
19& &   &   &    &    &   & 5\\
20& &   &   &    &    &   & 2\\
\hline
Total & 1 & 5 & 18 & 124 & 866 & 7622 & 72256\\
\hline
\end{tabular}

\end{center}

\address{Department of Mathematics \\
University of \AA rhus \\
8000 \AA rhus C \\
Denmark
}
{oebro@imf.au.dk}

\end{document}